\documentclass[a4paper,12pt]{article}

\usepackage{amsmath,amssymb}

\begin{document}

\begin{center}

\textbf{\Large Hypergroups arising from characters of } 

\medskip

\textbf{\Large  a compact group and its subgroup}

\medskip

\vspace{0,2 cm}

{\large Herbert Heyer, Satoshi Kawakami, Tatsuya Tsurii }

\medskip

{\large and Satoe Yamanaka}

\medskip


\end{center}

\begin{abstract}

The purpose of the present paper is to investigate a hypergroup 
arising from irreducible characters of a compact group $G$ 
and a closed subgroup $G_0$ of $G$ with index $[G:G_0]< + \infty$.  The convolution of 
this hypergroup is introduced by inducing irreducible representations 
of $G_0$ to $G$ and by restricting irreducible representations of 
$G$ to $G_0$. The method of proof relies on  
character formulae of  
induced representations of compact groups and of Frobenius' reciprocity theorem.

\bigskip

\medskip



\end{abstract}

\bigskip

\textbf{\large 1. Introduction}

\bigskip

One of the most challenging problems in the theory of hypergroups is a 
definite explanation of their algebraic structure. To solve 
this problem completely might be an utopian undertaking. 
But there are various ways to tackle parts of the problem. 
The approaches available are based on constructing new hypergroups from 
known ones. Much work has been done in the direction of extending 
hypergroups and of establishing new hypergroup structures defined by 
hypergroup actions ([HK1]).  
In succession of the authors' publications on semi-direct product 
hypergroups ([HK2]) and on hypergroup structures arising 
from certain dual objects of a hypergroup ([HK3]).   
The next step taken in the present paper is the supply of a 
hypergroup structure on the set of irreducible characters of a 
compact group $G$ together with a closed subgroup $G_0$ of $G$. 
It turns out that the resulting hypergroup structure can be 
characterized in terms of an invariance condition on characters 
of irreducible representations of $G_0$. 

\medskip

The method chosen in order to establish this result depends on 
the application of a character formula ([H]), of Frobenius' reciprocity 
theorem ([F]) for compact groups, and on the recently developed 
character theory for induced representations of hypergroups ([HKY]). 

\medskip

It should be mentioned that further progress in the 
research on the structure of hypergroups is on its way to publication : 
an extension of the notion of hyperfields ([HKKK]) to 
not necessarily finite hypergroups ([HKTY1]), and a 
generalization of the present work  to compact {\it hypergroups} and 
their closed {\it subhypergroups} ([HKTY2]). 

\medskip

A brief layout of the paper seems to be in order. 

\medskip

In Section 2 the preliminaries are restricted  
to the main notions of hypergroup theory ([BH], [J]) ; 
they can also be picked up from the introduction of [HK1]. 
In the  latter reference semi-direct product hypergroups have been 
introduced. 

\medskip

Section 3 is devoted to defining admissible pairs ($G,G_0$)  formed  
by a second countable compact group $G$ and a closed subgroup 
$G_0$ of $G$ of finite index, and to studying properties of 
these pairs (Lemma 3.5, 3.7, 3.9 and 3.10). 
In the Theorem of the section the set $\mathcal{K}(\hat{G} \cup \widehat{G_0})$ 
of characters of the union $\hat{G} \cup \widehat{G_0}$ of the 
duals of $G$ and $G_0$ is characterized as a hypergroup by 
the fact that ($G,G_0$) is an admissible pair. 
Applications of the Theorem are given to symmetric groups (Corollary 3.14) 
and to semi-direct product groups (Corollary 3.17). 

\medskip

Section 4 contains a variety of examples comprising those given in [SW].  
An extended list of examples can be visualized by Frobenius diagrams, 
a graph-theoretical illustration related to Dynkin diagrams 
and Coxeter graphs ([GHJ]). 
In all examples (except Example 4.6) the convolution of the hypergroup 
$\mathcal{K}(\hat{G} \cup \widehat{G_0})$ is 
described explicitly.    

\bigskip

\textbf{\large 2. Preliminaries} 

\medskip

For a locally compact space $X$ we shall mainly consider the subspaces $C_c(X)$ and $C_0(X)$ of the space $C(X)$ 
of  continuous functions on $X$ which have compact support or vanish at infinity respectively. 
By $M(X)$, $M^b(X)$ and $M_c(X)$ we abbreviate the spaces of all (Radon) measures on $X$, the 
bounded measures and the measures with compact support on $X$ respectively. 
Let $M^1(X)$ denote the set of probability measures on $X$ and $M^1_c(X)$ its subset $M^1(X) \cap M_c(X)$.
The symbol $\delta_x$ stands for the Dirac measures in $x \in X$. 

\medskip

A {\it hypergroup} $(K,*)$ is a locally compact space $K$ together with a {\it convolution} $*$ in $M^b(K)$ such 
that $(M^b(K),*)$ becomes a Banach algebra and that the following properties are fulfilled. 

\medskip

(H1) The mapping 
$$
(\mu, \nu )\longmapsto \mu *  \nu
$$

\ \ \ \ \ \ \ from $M^b(K) \times M^b(K)$ into $M^b(K)$ is continuous with respect to the 

\ \ \ \ \ \ \ weak  topology in $M^b(K)$. 

\medskip

(H2) For $x,y \in K$ the convolution $\delta_x * \delta_y$ belongs to $M^1_c(K)$.

\medskip

(H3) There exists a {\it unit} element $e \in K$ with 
$$
\delta_e * \delta_x = \delta_x * \delta_e = \delta_x 
$$

\ \ \ \ \ \ \ for all $x \in K$, and an {\it involution} 
$$
x \longmapsto x^-
$$
   
\ \ \ \ \ \ \ in $K$ such that 
$$
\delta_{x^-} * \delta_{y^-} = (\delta_y * \delta_x)^-
$$

\ \ \ \ \ \ \ and 
$$
e \in \text{supp}(\delta_x * \delta_y)~~\text{if and only if}~~x = y^-
$$

\ \ \ \ \ \ \ whenever $x,y \in K$. 

\medskip

(H4) The mapping 
$$
(x,y) \longmapsto \text{supp}(\delta_x * \delta_y)
$$

\ \ \ \ \ \ \ from $K \times K$ into the space $\mathcal{C}(K)$ of all compact subsets of $K$ furnished  

\ \ \ \ \ \ \ with Michael topology is continuous. 

\medskip

A hypergroup $(K, *)$ is said to be {\it commutative} if the convolution $*$ is commutative. 
In this case $(M^b(K), *, -)$ is a commutative Banach $*$-algebra with identity $\delta_e$. 
There is an abundance of hypergroups and there are various constructions (polynomial, Sturm-Liouville) as the reader 
may learn from the pioneering papers on the subject and also from the monograph [BH]. 

\medskip

Let $(K,*)$ and $(L, \circ)$ be two hypergroups with units $e_K$ and $e_L$ respectively. A continuous mapping 
$\varphi:K \rightarrow L$ is called a hypergroup {\it homomorphism} if $\varphi(e_K)=  e_L$ and $\varphi$ is the unique linear, weakly continuous 
extension from $M^b(K)$ to $M^b(L)$ 
such that 
$$
\varphi(\delta_x) = \delta_{\varphi(x)}, \ 
\varphi(\delta_x^-) = \varphi(\delta_x)^- \  \textrm{and} \  
\varphi (\delta_x * \delta_y) = \varphi (\delta_x) \circ \varphi (\delta_y)
$$
whenever $x,y \in K$. 
If $\varphi:K \rightarrow L$ is also a homeomorphism, it will be called 
an {\it isomorphism} from $K$ onto $L$. 
An isomorphism from $K$ onto $K$ is called an {\it automorphism} of $K$. 
We denote by Aut($K$) the set of all automorphisms of $K$. 
Then Aut($K$) becomes a topological group equipped with the weak topology 
of $M^b(K)$. 
We call $\alpha$ an \textit{action} of a locally compact group $G$ on a hypergroup 
$H$ if $\alpha$ is a continuous homomorphism from $G$ into Aut($H$). 
Associated with the action $\alpha$ of $G$ on $H$ one can define a 
semi-direct product hypergroup $K=H\rtimes_{\alpha}G$, see [HK1].  

\medskip

\medskip

If the given hypergroup $K$ is commutative, its dual $\widehat{K}$ can be  introduced as the set of all bounded continuous functions 
$\chi \not = 0$ on $K$ satisfying
$$
\int_K \chi(z)(\delta_x^- * \delta_y)(dz) = \overline{\chi(x)}\chi(y)
$$
for all $x,y \in K$. This set of characters $\widehat{K}$ of $K$ 
becomes a locally compact space with  respect to 
the topology of uniform convergence on compact sets, 
but generally fails to be a hypergroup. 
When $\widehat{K}$ is 
a hypergroup $K$ is called a strong hypergroup. When the dual 
$\widehat{K}$ of a strong hypergroup $K$ is also strong and 
$\widehat{\widehat{K}} \cong K$ holds  $K$ is called a Pontryagin hypergroup. 

\medskip

\bigskip

\textbf{\large 3. Hypergroups related to admissible pairs}

\medskip

Let $G$ be a compact group which satisfies the second axiom of 
countability and let $\hat{G}$ be the set of all equivalence classes of 
irreducible representations of $G$. Then $\hat{G}$ is (at most  countable)  
discrete space which we write explicitly as  
$$
\hat{G} = \{\pi_0, \pi_1, \cdots, \pi_n,\cdots \}, 
$$
where $\pi_0$ is the trivial representation of $G$. 
We denote by ${\rm Rep}^{\rm f}(G)$ the set of equivalence classes of 
finite-dimensional representations of $G$. For  $\pi \in {\rm Rep}^{\rm f}(G)$  
we consider the normalized character of $\pi$ given by 
$$
ch(\pi)(g) = \frac{1}{\dim \pi}{\rm tr}(\pi(g)),  
$$
for all $g \in G$.

Put 
$$
\mathcal{K}(\hat{G}) = \{ch(\pi) : \pi \in \hat{G}\}. 
$$
Then $\mathcal{K}(\hat{G})$ is known to be a discrete commutative hypergroup with 
unit $ch(\pi_0) = \pi_0$.

\medskip

Let $G_0$ be a closed subgroup of $G$ such that the index $[G:G_0] $ is finite. 
We write $\widehat{G_0} = \{\tau_0, \tau_1, \cdots, \tau_n, \cdots\}$,  
where $\tau_0$ is the trivial representation of $G_0$.

\medskip

The following is well-known fact. 

\medskip
\noindent
{\bf Lemma 3.1  } 

\medskip

(1) $ch(\pi_i \otimes \pi_j) = ch(\pi_i) ch(\pi_j)$ 
for $\pi_i, \pi_j \in  {\rm Rep}^{\rm f}(G)$. 

\medskip

(2) $ch({\rm res}_{G_0}^G \pi) = {\rm res}_{G_0}^G ch(\pi)$ 
for $\pi \in {\rm Rep}^{\rm f}(G)$.

\medskip
\noindent
{\bf Lemma 3.2 }  

\medskip

(1) \ [Character formula ] (see Hirai [H])\ For $\tau \in {\rm Rep}^{\rm f}(G_0)$, 
$$
ch({\rm ind}_{G_0}^G \tau)(g) = \int_G ch(\tau)(sgs^{-1})  
1_{G_0}(sgs^{-1}) \omega_G(ds). 
$$

(2) \ [Frobenius' reciprocity theorem] (see Folland [F]) \ For $\tau \in \widehat{G_0}$ 
and   for $\pi \in \hat{G}$,  
$$ 
[{\rm ind}_{G_0}^G \tau : \pi ] = [\tau : {\rm res}_{G_0}^G \pi],  
$$
where $[ \ \ : \ \ ]$ denotes the multiplicity of representations. 

\medskip
\noindent
{\bf Remark } $ch({\rm ind}_{G_0}^G \tau) = {\rm ind}_{G_0}^G ch(\tau)$, see [HKY], where 
$$
{\rm ind}_{G_0}^G ch(\tau) := \int_G ch(\tau)(sgs^{-1})  
1_{G_0}(sgs^{-1}) \omega_G(ds).
$$

\medskip
\noindent
{\bf Definition } On the set 
$$
\mathcal{K}(\hat{G} \cup \widehat{G_0})
:= \{(ch(\pi), \circ), (ch(\tau), \bullet) : \pi \in \hat{G}, \tau \in \widehat{G_0}\}  
$$
we define a convolution $*$ as follows. 
For $\pi_i, \pi_j, \pi \in \hat{G}$ and 
$\tau_i, \tau_j, \tau \in \widehat{G_0}$, 
\begin{align*}
&(ch(\pi_i), \circ) * (ch(\pi_j), \circ) := (ch(\pi_i)ch(\pi_j), \circ),\\
&(ch(\pi), \circ) * (ch(\tau), \bullet) := (ch({\rm res}_{G_0}^G \pi) ch(\tau), \bullet),\\
&(ch(\tau), \bullet) * (ch(\pi), \circ) := (ch(\tau)ch({\rm res}_{G_0}^G \pi), \bullet),\\
&(ch(\tau_i), \bullet) * (ch(\tau_j), \bullet) 
:= (ch({\rm ind}_{G_0}^G(\tau_i \otimes \tau_j)), \circ).
\end{align*}
We want to check the associativity relations of the convolution in the following cases. 
Whenever reference to a particular representation $\pi$ is not needed, 
we abbreviate  $(ch(\pi), \bullet)$ by $\bullet$ and $(ch(\pi), \circ)$ 
by $\circ$. Hence our task will be to verify the subsequent formulae : 
\begin{align*}
&(A1)~(\circ * \circ) * \circ = \circ * (\circ * \circ), \\
&(A2)~(\bullet * \circ) * \circ = \bullet * (\circ * \circ),\\
&(A3)~(\bullet * \bullet) * \circ = \bullet * (\bullet * \circ)~{\rm and}\\
&(A4)~(\bullet * \bullet) * \bullet = \bullet * (\bullet * \bullet). 
\end{align*}

\medskip
\noindent
{\bf Lemma 3.3 } The equalities $(A1)$, $(A2)$ and $(A3)$ hold 
without further assumptions.  
For $\pi_i, \pi_j, \pi_k, \pi \in \hat{G}$ and 
$\tau_i, \tau_j, \tau \in \widehat{G_0}$, 
\begin{align*}
&(A1)~((ch(\pi_i), \circ) * (ch(\pi_j), \circ)) * (ch(\pi_k), \circ) 
= (ch(\pi_i), \circ) * ((ch(\pi_j), \circ) * (ch(\pi_k), \circ)). \\
&(A2)~ ((ch(\tau), \bullet) * (ch(\pi_i), \circ)) * (ch(\pi_j), \circ)
=  (ch(\tau), \bullet) * ((ch(\pi_i), \circ) * (ch(\pi_j), \circ)).\\
&(A3)~((ch(\tau_i), \bullet) * (ch(\tau_j), \bullet)) * (ch(\pi), \circ)
=(ch(\tau_i), \bullet) * ((ch(\tau_j), \bullet) * (ch(\pi), \circ)).
\end{align*}

\medskip
\noindent
{\bf Proof } $(A1)$ is clear because $\mathcal{K}(\hat{G})$ has a hypergroup structure. 

\medskip

$(A2)$ For $\tau \in \widehat{G_0}$ and $\pi_i, \pi_j \in \hat{G}$, 
\begin{align*}
&((ch(\tau), \bullet) * (ch(\pi_i), \circ)) * (ch(\pi_j), \circ)\\
=& (ch(\tau) ch({\rm res}_{G_0}^G \pi_i), \bullet) * (ch(\pi_j), \circ)\\
=& (ch(\tau) ch({\rm res}_{G_0}^G \pi_i) ch({\rm res}_{G_0}^G \pi_j) , \bullet)\\
=& (ch(\tau) ({\rm res}_{G_0}^G ch(\pi_i)) ({\rm res}_{G_0}^G ch(\pi_j)) , \bullet). 
\end{align*}
On the other hand, 
\begin{align*}
&(ch(\tau), \bullet) * ((ch(\pi_i), \circ) * (ch(\pi_j), \circ))\\
=& (ch(\tau), \bullet)) * ((ch(\pi_i) ch(\pi_j), \circ)\\
=& (ch(\tau) {\rm res}_{G_0}^G (ch(\pi_i) ch(\pi_j)), \bullet) \\
=& (ch(\tau) ({\rm res}_{G_0}^G ch(\pi_i)) ({\rm res}_{G_0}^G ch(\pi_j)), \bullet).
\end{align*}

\medskip
(A3) For $\tau_i, \tau_j \in \widehat{G_0}$ and $\pi \in \hat{G}$, 
\begin{align*}
&((ch(\tau_i), \bullet) * (ch(\tau_j), \bullet)) * (ch(\pi), \circ)\\
=&
 (ch({\rm ind}_{G_0}^G(\tau_i \otimes \tau_j)), \circ) * (ch(\pi), \circ)\\
=& (ch({\rm ind}_{G_0}^G(\tau_i \otimes \tau_j)) ch(\pi), \circ). 
\end{align*}
For every $g \in G$, 
\begin{align*}
&(ch({\rm ind}_{G_0}^G(\tau_i \otimes \tau_j)) ch(\pi))(g) \\
=&  ch({\rm ind}_{G_0}^G(\tau_i \otimes \tau_j)(g) ch(\pi) (g)\\
= & \int_G ch(\tau_i \otimes \tau_j)(sgs^{-1}) 
1_{G_0}(sgs^{-1}) \omega_G(ds)  ch(\pi) (g) \\
= & \int_G (ch(\tau_i) ch(\tau_j))(sgs^{-1}) ch(\pi)(g)
1_{G_0}(sgs^{-1}) \omega_G(ds)  \\
= & \int_G ch(\tau_i)(sgs^{-1}) ch(\tau_j)(sgs^{-1})  ch(\pi) (sgs^{-1})
1_{G_0}(sgs^{-1})\omega_G(ds). 
\end{align*}
On the other hand, 
\begin{align*}
&(ch(\tau_i), \bullet) * ((ch(\tau_j), \bullet) * (ch(\pi), \circ))\\
=& (ch(\tau_i), \bullet) * (ch(\tau_j)ch({\rm res}_{G_0}^G \pi), \bullet)\\
=& (ch({\rm ind}_{G_0}^G (\tau_i \otimes \tau_j \otimes {\rm res}_{G_0}^G \pi)), \circ).   
\end{align*}
For each $g \in G$, 
\begin{align*}
&ch({\rm ind}_{G_0}^G (\tau_i \otimes \tau_j \otimes {\rm res}_{G_0}^G \pi))(g) \\
= & \int_G ch(\tau_i \otimes \tau_j \otimes {\rm res}_{G_0}^G \pi)(sgs^{-1}) 
1_{G_0}(sgs^{-1}) \omega_G(ds)\\
= & \int_G (ch(\tau_i)ch(\tau_j)ch({\rm res}_{G_0}^G \pi))(sgs^{-1}) 
1_{G_0}(sgs^{-1}) \omega_G(ds)\\
= & \int_G ch(\tau_i)(sgs^{-1})ch(\tau_j)(sgs^{-1})ch(\pi)(sgs^{-1}) 
1_{G_0}(sgs^{-1}) \omega_G(ds). 
\end{align*}
\hspace{\fill}[Q.E.D.]

\medskip
\noindent
{\bf Definition } Let $(G, G_0)$ be a pair of consisting of a compact group $G$ and 
a closed subgroup $G_0$ of $G$. 
For $g \in G_0$
$$
X(g):= \{ s\in G : sgs^{-1} \in G_0 \}.
$$ 
We call $(G, G_0)$ an  {\it  admissible pair}   
if for any $\tau \in \widehat{G_0}$, any $g \in G_0$ and any $s \in X(g)$, 
$$
ch(\tau)(sgs^{-1}) = ch(\tau)(g)
$$
holds.  

\medskip
\noindent
{\bf Lemma 3.4 } If a compact group $G$ together with a subgroup $G_0$ 
of $G$ with $[G:G_0] < + \infty$ forms an admissible pair, then the associativity relation $(A4)$ holds.

\medskip
\noindent
{\bf Proof } Assume that $(G, G_0)$ is an admissible pair. For $\tau_i, \tau_j, \tau_k \in \widehat{G_0}$ 
and $g \in G_0$ 
\begin{align*}
&((ch(\tau_i), \bullet) * (ch(\tau_j), \bullet)) * (ch(\tau_k), \bullet)\\
=&(ch({\rm ind}_{G_0}^G (\tau_i \otimes \tau_j)), \circ) * (ch(\tau_k), \bullet)\\
=&({\rm ind}_{G_0}^G ch(\tau_i \otimes \tau_j)), \circ) * (ch(\tau_k), \bullet)\\
=&({\rm res}_{G_0}^G({\rm ind}_{G_0}^G (ch(\tau_i) ch(\tau_j))) (ch(\tau_k), \bullet). 
\end{align*}
For $g \in G_0$, 
\begin{align*}
&({\rm ind}_{G_0}^G (ch(\tau_i) ch(\tau_j)))(g) ch(\tau_k))(g)\\
=& \left(\int_{G} ch(\tau_i)(sgs^{-1}) ch(\tau_j)(sgs^{-1}) 1_{G_0}(sgs^{-1})  
\omega_G(ds)\right) ch(\tau_k)(g)\\
=& \left(\int_{G} ch(\tau_i)(g) ch(\tau_j)(g) 1_{G_0}(sgs^{-1})  
\omega_G(ds)\right) ch(\tau_k)(g)\\
=& \left(\int_{G} 1_{G_0}(sgs^{-1})  
\omega_G(ds)\right) ch(\tau_i)(g) ch(\tau_j)(g) ch(\tau_k)(g). 
\end{align*}
This implies the associativity relation $(A4)$. 
\hspace{\fill}[Q.E.D.]

\medskip
\noindent
{\bf Lemma 3.5 } If  the associativity relation $(A4)$ holds  
for a compact group $G$ and a subgroup $G_0$ of $G$ with $[G:G_0] < + \infty$, 
then $(G,G_0)$  is an admissible pair. 

\medskip
\noindent
{\bf Proof }
Assume that the associativity relation $(A4)$ holds. Let $\tau_0$ be the trivial 
representation of $\widehat{G_0}$. For $\tau \in \widehat{G_0}$ the associativity 
relation  
$$
((ch(\tau_0), \bullet) * (ch(\tau_0), \bullet)) * (ch(\tau), \bullet) 
= (ch(\tau_0), \bullet) * ((ch(\tau_0), \bullet) * (ch(\tau), \bullet)) 
$$
holds. 
\begin{align*}
((ch(\tau_0), \bullet) * (ch(\tau_0), \bullet)) * (ch(\tau), \bullet)
= (ch({\rm res}_{G_0}^G ({\rm ind}_{G_0}^G  \tau_0 )) ch(\tau),\bullet) 
\end{align*}
and
$$
(ch(\tau_0), \bullet) * ((ch(\tau_0), \bullet) * (ch(\tau), \bullet))
= (ch({\rm res}_{G_0}^G ({\rm ind}_{G_0}^G \tau)), \bullet). 
$$
Then for $g \in G_0$
$$
ch({\rm ind}_{G_0}^G \tau_0)(g) ch(\tau)(g) = ch({\rm ind}_{G_0}^G \tau)(g). 
$$
Now  
$$
ch({\rm ind}_{G_0}^{G} \tau_0)(g) \geq \omega_{G}(G_0) >0. 
$$
Indeed by the character formula 
$$
ch({\rm ind}_{G_0}^{G} \tau_0)(g)  
= \int_{G}  1_{G_0} (sgs^{-1})  \omega_{G}(ds) = \omega_{G}(X(g)). 
$$
Since $X(g) \supset G_0$, we see that $\omega_{G}(X(g)) \geq \omega_{G}(G_0)$. 
By the assumption  $[G :G_0] < + \infty $ we obtain $\omega_{G}(G_0) = 1/[G :G_0] > 0$. 
Hence 
$$
ch(\tau)(g) = (ch({\rm ind}_{G_0}^G \tau_0)(g))^{-1} ch({\rm ind}_{G_0}^G \tau)(g). 
$$
For $s \in X(g)$
\begin{align*}
ch(\tau)(sgs^{-1}) 
&= (ch({\rm ind}_{G_0}^{G} \tau_0)(sgs^{-1}))^{-1}  ch({\rm ind}_{G_0}^{G} \tau)(sgs^{-1})\\
&= (ch({\rm ind}_{G_0}^{G} \tau_0)(g))^{-1}  ch({\rm ind}_{G_0}^{G} \tau)(g)\\
&= ch(\tau)(g).  
\end{align*}
Then $(G, G_0)$ is an admissible pair. 
\hspace{\fill}[Q.E.D.]

\bigskip
\noindent
{\bf Theorem } Let $G_0$ be a closed subgroup of a compact group $G$ such that 
$[G:G_0] < + \infty$. Then $\mathcal{K}(\hat{G} \cup \widehat{G_0})$ is a hypergroup 
if and only if $(G, G_0)$ is an admissible pair.

\medskip
\noindent
{\bf Proof } The associativity relations $(A1)$, $(A2)$ and $(A3)$ are 
a consequence of Lemma 3.3, and $(A4)$ holds if and only if 
$(G, G_0)$ is an admissible pair by Lemma 3.4 and  Lemma 3.5. 
It is easy to check the remaining axioms of a hypergroup 
for $\mathcal{K}(\hat{G} \cup \widehat{G_0})$. 
The desired conclusion follows. 
\hspace{\fill}[Q.E.D.]

\medskip
\noindent
{\bf Remark 3.6 } 

\medskip
(1) The above $\mathcal{K}(\hat{G} \cup \widehat{G_0})$ is a discrete commutative 
(at most countable) hypergroup such that the sequence : 
$$
1 \longrightarrow \mathcal{K}(\hat{G}) \longrightarrow \mathcal{K}(\hat{G} \cup \widehat{G_0}) 
\longrightarrow \mathbb{Z}_2 \longrightarrow 1
$$
is exact. 

\medskip

(2) If $G_0 = G$, then $\mathcal{K}(\hat{G} \cup \widehat{G_0})$ is the hypergroup 
$\mathcal{K}(\hat{G}) \times \mathbb{Z}_2$. 

\medskip

(3) If $G$ is a finite group and $G_0 = \{e\}$, then  
$\mathcal{K}(\hat{G} \cup \widehat{G_0})$ is the hypergroup join 
$\mathcal{K}(\hat{G}) \vee \mathbb{Z}_2$.

\medskip
\noindent
{\bf Lemma 3.7 } If $G$ is a compact Abelian group and $G_0$ a closed subgroup of $G$ 
with $[G: G_0] < + \infty$. Then $(G, G_0)$ is always an admissible pair. 

\medskip
\noindent
{\bf Proof } The desired assertion clearly follows from the fact that $sgs^{-1} = g$ 
for $g \in G_0$ and $s \in G$. 
\hspace{\fill}[Q.E.D.]

\bigskip
\noindent
{\bf Corollary 3.8 } Let $G$ be a compact Abelian group and $G_0$ a closed subgroup of $G$ 
with $[G: G_0] < + \infty$. Then $\mathcal{K}(\hat{G} \cup \widehat{G_0})$ 
is a hypergroup. 

\medskip
\noindent
{\bf Proof } This assertion follows directly from the Theorem  and Lemma 3.7. 
\hspace{\fill}[Q.E.D.]

\medskip
\noindent
{\bf Lemma 3.9 } If for each $\tau \in \widehat{G_0}$ there exists a representation 
$\tilde{\tau}$ of $G$ such that ${\rm res}_{G_0}^G \tilde{\tau} = \tau$, then 
$(G, G_0)$ is an admissible pair.

\medskip
\noindent
{\bf Proof }  For $\tau \in \widehat{G_0}$, $g \in G_0$ and $s \in X(g)$, 
$$
ch(\tau)(sgs^{-1}) = ch(\tilde{\tau})(sgs^{-1}) = ch(\tilde{\tau})(g) = ch(\tau)(g). 
$$
\hspace{\fill}[Q.E.D.]

\bigskip
\noindent
{\bf Corollary 3.10 } Let $G$ be a semi-direct product group 
$H \rtimes_\alpha G_0$,  where $H$ is a finite group and $G_0$ is a 
finite group. Then 
$\mathcal{K}(\hat{G} \cup \widehat{G_0})$ is a hypergroup.

\medskip
\noindent
{\bf Proof } For $\tau \in \widehat{G_0}$, put 
$$
\tilde{\tau}((h,g)) = \tau(g)
$$
for $(h,g) \in H \rtimes_\alpha G_0 = G$. Then $\tilde{\tau}$ is a finite 
dimensional representation of $G$ and ${\rm res}_{G_0}^G \tilde{\tau} = \tau.$ 
By the Theorem  and Lemma 3.9 we arrive at the desired conclusion. 

\hspace{\fill}[Q.E.D.]

\medskip
\noindent
{\bf Lemma 3.11 } If for $g \in G_0$ and $s \in X(g)$ there exists $t \in G_0$ 
such that $tgt^{-1} = sgs^{-1}$, then $(G, G_0)$ is an admissible pair. 

\medskip
\noindent
{\bf Proof } For $ \tau \in \widehat{G_0}$, $g \in G_0$ and $s \in X(g)$, 
$$
ch(\tau)(sgs^{-1}) = ch(\tau)(tgt^{-1}) = ch(\tau)(g). 
$$  
\hspace{\fill}[Q.E.D.]

\bigskip

Let $S_n$ be the symmetric group of degree $n$.

\medskip

\noindent
{\bf Corollary 3.12 } $\mathcal{K}(\widehat{S_n} \cup \widehat{S_{n-1}})$ $(n \geq 2)$ is a 
hypergroup.

\medskip
\noindent
{\bf Proof }  For  $g \in S_{n-1}$, $s \in X(g)$ such that $s^{-1}(n) = a$, 
$$
sgs^{-1}(n) = s(g(s^{-1}(n))) = s(g(a)). 
$$
Since $sgs^{-1} \in S_{n-1}$, $sgs^{-1}(n) = n$ and $s(g(a)) = n$. 
Then $g(a) = s^{-1}(n) = a$. 

Put $t = s s_1$ where $s_1$ is a transposed permutation $(a, n)$. Then 
we  see that 
$$
tgt^{-1} = sgs^{-1} 
$$
by the fact that for $b$ such that $b \not = a$, $s^{-1}(b) \not = a$ and 
$g(s^{-1}(b)) \not = a$ hold. By the Theorem  and Lemma 3.11 we get the desired conclusion. 
\hspace{\fill}[Q.E.D.]

\medskip
\noindent
{\bf Lemma 3.13 } Let $G_0$ and $G_1$ be  closed subgroups of $G$ such that 
$G_0 \subset G_1 \subset G$. If $(G_1, G_0)$ and $(G, G_1)$ are admissible pairs, 
then $(G, G_0)$ is an admissible pair.

\medskip
\noindent
{\bf Proof } Since $(G_1, G_0)$ is an admissible pair, for $\tau \in \widehat{G_0}$ and 
$g \in G_0$, 
\begin{align*}
ch({\rm ind}_{G_0}^{G_1} \tau)(g)  
&= \int_{G_1} ch(\tau)(sgs^{-1}) 1_{G_0} (sgs^{-1}) d \omega_{G_1}(s)\\
&= \int_{G_1} ch(\tau)(g) 1_{G_0} (sgs^{-1}) d \omega_{G_1}(s)\\
&= \int_{G_1}  1_{G_0} (sgs^{-1}) d \omega_{G_1}(s) ch(\tau)(g)\\
&= ch({\rm ind}_{G_0}^{G_1} \tau_0)(g) ch(\tau)(g), 
\end{align*}
where $\tau_0$ is the trivial representation of $G_0$. 
Since 
$$
ch({\rm res}_{G_0}^{G_1} \tau_0)(g) \geq \omega_{G_1}(G_0) >0,
$$  
we see that
$$
ch(\tau)(g) =  (ch({\rm ind}_{G_0}^{G_1} \tau_0)(g))^{-1}  ch({\rm ind}_{G_0}^{G_1} \tau)(g). 
$$
Since  $(G, G_1)$ is an admissible pair and ${\rm ind}_{G_0}^{G_1} \tau_0 \in {\rm Rep^f}(G_1)$,   
${\rm ind}_{G_0}^{G_1} \tau \in {\rm Rep^f}(G_1)$, we have for $g \in G_0 \subset G_1$ and 
$s \in X(g)$, 
$$
ch({\rm ind}_{G_0}^{G_1} \tau_0)(sgs^{-1}) = ch({\rm ind}_{G_0}^{G_1} \tau_0)(g)
$$
and 
$$
ch({\rm ind}_{G_0}^{G_1} \tau)(sgs^{-1}) = ch({\rm ind}_{G_0}^{G_1} \tau)(g). 
$$
Then we obtain 
\begin{align*}
ch(\tau)(sgs^{-1}) 
&= (ch({\rm ind}_{G_0}^{G_1} \tau_0)(sgs^{-1}))^{-1}  ch({\rm ind}_{G_0}^{G} \tau)(sgs^{-1})\\
&= (ch({\rm ind}_{G_0}^{G_1} \tau_0)(g))^{-1}  ch({\rm ind}_{G_0}^{G_1} \tau)(g)\\
&= ch(\tau)(g) 
\end{align*}
\hspace{\fill}[Q.E.D.]

\bigskip
\noindent
{\bf Corollary 3.14 } For natural numbers $m$ and $n$ such that 
$m > n \geq 1$,  $\mathcal{K}(\widehat{S_m} \cup \widehat{S_n})$ is a 
hypergroup.

\medskip
\noindent
{\bf Proof } This statement follows from the Theorem  and Lemma 3.13 combined with 
Corollary 3.12.  
\hspace{\fill}[Q.E.D.]

\bigskip
Let $G_0$ be a closed normal subgroup of $G$. Then the coadjoint action 
$\hat{\alpha}$ of $G$ on $\widehat{G_0}$ is defined by 
$$
\hat{\alpha}_s (\tau)(g) := \tau(sgs^{-1})
$$
for $\tau \in \widehat{G_0}$, $g \in G_0$ and $s \in G$. 
If $\hat{\alpha}_s = id $ for all $s \in G$, we say that  $\hat{\alpha}$ is trivial.

\medskip
\noindent
{\bf Lemma 3.15 } Let $G_0$ be a closed normal subgroup of $G$. The 
pair $(G, G_0)$ is an admissible pair if and only if the coadjoint action 
$\hat{\alpha}$ is trivial. 

\medskip
\noindent
{\bf Proof  } Assume that $(G, G_0)$ is an admissible pair. For $g \in G_0$ it 
is clear that $X(g) = G$. Then for $\tau \in \widehat{G_0}$
$$
ch(\tau)(sgs^{-1}) = ch(\tau)(g) 
$$ 
for all $s \in G$. This implies that 
$$
ch(\hat{\alpha}_s(\tau))(g) = ch(\tau)(g) 
$$
for all $g \in G_0$. Hence we obtain 
$$
\hat{\alpha}_s(\tau) \cong \tau
$$
for $\tau \in \widehat{G_0}$. In fact $\hat{\alpha}_s$ is the identity on $\widehat{G_0}$ which 
means that $\hat{\alpha}$ is trivial. 

The converse is clear.  
\hspace{\fill}[Q.E.D.]

\bigskip
\noindent
{\bf Lemma 3.16 } Let $G_0$ be a closed normal commutative subgroup of $G$. 
The pair $(G, G_0)$ is an admissible pair if and only if $G \cong G_0 \times (G/{G_0})$.

\bigskip
\noindent
{\bf Proof  } Assume that $(G, G_0)$ is an admissible pair and $sgs^{-1} \not = g$ for 
$g \in G_0$ and $s \in X(g) = G$. Since $\widehat{G_0}$ separates $G_0$, there exists 
$\tau \in \widehat{G_0}$ such that
$$
\tau(sgs^{-1}) \not = \tau(g). 
$$ 
This contradicts the assumption that $(G, G_0)$ is an admissible pair. 
But then for $g \in G_0$ and $s \in G$
$$
sgs^{-1} = g, 
$$
namely
$$
sg = gs
$$
holds. For $g_1, g_2 \in G_0$ and $s_1, s_2 \in G$
$$
(g_1 s_1) (g_2 s_2) = g_1 (s_1 g_2 ) s_2 = g_1(g_2 s_1) s_2 = (g_1 g_2)(s_1 s_2).
$$
This implies that $G \cong G_0 \times (G/G_0)$. 

The converse is clear by Lemma 3.15. 
\hspace{\fill}[Q.E.D.]

\bigskip
\noindent
{\bf Corollary 3.17 } Let $G$ be a semi-direct product group $H \rtimes_\alpha G_0$ 
where  $H$ is a compact Abelian group and $G_0$ is a finite group. 
$\mathcal{K}(\hat{G} \cup \hat{H})$ is a hypergroup if and only if 
the action $\alpha$ is trivial, i.e. $G = H \times G_0$. 

\medskip
\noindent
{\bf Proof } We note that $H$ is a closed normal subgroup of $G$ and $G/H \cong G_0$. 
Then the assertion follows from the Theorem together with Lemma 3.16. 
\hspace{\fill}[Q.E.D.]

\bigskip

\textbf{\large 4. Examples }

\medskip

Associated with a pair $(G, G_0)$ of finite groups such that $G \supset G_0$, 
we obtain a certain finite graph $D(\hat{G} \cup \widehat{G_0})$ by Frobenius' 
reciprocity theorem. The set of vertices is  
$$
 \{(\pi, \circ), (\tau, \bullet) : \pi \in \hat{G}, \tau \in \widehat{G_0}\}  
$$
and the edge between $(\pi, \circ)$ and $(\tau, \bullet)$ is given by the 
multiplicity 
$$
m_{\pi,\tau} := [{\rm ind}_{G_0}^G(\tau) : \pi ] = [\tau : {\rm res}_{G_0}^G \pi] \neq 0.
$$
We call this graph $D(\hat{G} \cup \widehat{G_0})$ a Frobenius diagram. 
Frobenius diagrams $D(\hat{G} \cup \widehat{G_0})$ sometimes appear 
as Dynkin diagrams and sometimes as Coxter graphs ([GHJ]). V. S. Sunder and N. J. Wildberger 
 constructed in [SW] fusion rule algebras $\mathcal{F}(D)$ and hypergroups 
 $\mathcal{K}(D)$ associated with certain Dynkin diagrams of type $A_n$,  $D_{2n}$ and so on. 
 We give some examples of $\mathcal{K}(\hat{G} \cup \widehat{G_0})$ 
which are compatible with Frobenius diagrams $D(\hat{G} \cup \widehat{G_0})$.

\bigskip
\noindent
{\bf 4.1 } \ The case that 
$G = \mathbb{Z}_2 = \{e,g\}$ $(g^2 = e)$ and $G_0 = \{e\}$. 
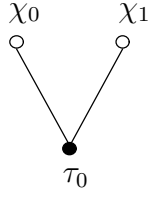
\begin{center}
{\unitlength 0.1in%
\begin{picture}(  6.2400,  8.5900)( 19.8000,-17.0000)%
%
\special{pn 8}%
\special{pa 2021 1109}%
\special{pa 2295 1611}%
\special{fp}%
\special{pa 2302 1604}%
\special{pa 2570 1109}%
\special{fp}%
\put(19.8000,-9.7100){\makebox(0,0)[lb]{$\chi_0$}}%
\put(25.4900,-9.7100){\makebox(0,0)[lb]{$\chi_1$}}%
\put(22.6100,-18.3000){\makebox(0,0)[lb]{$\tau_0$}}%
%
\special{pn 8}%
\special{ar 2570 1072 34 32  0.0000000  6.2831853}%
%
\special{pn 8}%
\special{ar 2021 1072 35 32  0.0000000  6.2831853}%
%
\special{sh 1.000}%
\special{ia 2295 1629 35 32  0.0000000  6.2831853}%
\special{pn 8}%
\special{ar 2295 1629 35 32  0.0000000  6.2831853}%
\end{picture}}%
\end{center}

\bigskip

$\mathcal{K}(\hat{G} \cup \widehat{G_0}) = \{(ch(\chi_0), \circ), (ch(\chi_1), \circ), 
(ch(\tau_0), \bullet)\}$. Put $\gamma_0 = (ch(\chi_0), \circ)$, 
$\gamma_1 = (ch(\chi_1), \circ)$ and $\rho_0 = (ch(\tau_0), \bullet)$. Then 
the structure equations are
\begin{align*}
\gamma_1 \gamma_1 = \gamma_0,~~~\rho_0 \rho_0 = \frac{1}{2}\gamma_0 + \frac{1}{2}\gamma_1,~~~
\gamma_1 \rho_0 = \rho_0. 
\end{align*}

\bigskip
\noindent
{\bf 4.2 } \ The case that 
$G = \mathbb{Z}_3 = \{e,g, g^2\}$ $(g^3 = e)$ and $G_0 = \{e\}$. 
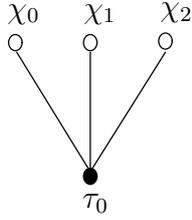
\begin{center}
{\unitlength 0.1in%
\begin{picture}(  8.5300,  9.9300)( 19.8000,-18.9000)%
\put(19.8000,-10.3500){\makebox(0,0)[lb]{$\chi_0$}}%
\put(23.7500,-10.3500){\makebox(0,0)[lb]{$\chi_1$}}%
\put(23.6800,-20.2000){\makebox(0,0)[lb]{$\tau_0$}}%
%
\special{pn 8}%
\special{ar 2409 1133 35 42  0.0000000  6.2831853}%
%
\special{pn 8}%
\special{ar 2021 1133 35 42  0.0000000  6.2831853}%
%
\special{sh 1.000}%
\special{ia 2409 1831 35 42  0.0000000  6.2831853}%
\special{pn 8}%
\special{ar 2409 1831 35 42  0.0000000  6.2831853}%
%
\special{pn 8}%
\special{ar 2798 1125 35 42  0.0000000  6.2831853}%
%
\special{pn 8}%
\special{pa 2409 1183}%
\special{pa 2409 1807}%
\special{fp}%
\special{pa 2798 1183}%
\special{pa 2409 1807}%
\special{fp}%
\special{pa 2028 1183}%
\special{pa 2409 1807}%
\special{fp}%
\put(27.7100,-10.2700){\makebox(0,0)[lb]{$\chi_2$}}%
\end{picture}}%
\end{center}

\bigskip

$\mathcal{K}(\hat{G} \cup \widehat{G_0}) = \{(ch(\chi_0), \circ), (ch(\chi_1), \circ), 
(ch(\chi_2), \circ), (ch(\tau_0), \bullet)\}$. Put $\gamma_0 = (ch(\chi_0), \circ)$, 
$\gamma_1 = (ch(\chi_1), \circ)$, $\gamma_2 = (ch(\chi_2), \circ)$ 
and $\rho_0 = (ch(\tau_0), \bullet)$. Then 
the structure equations are
\begin{align*}
&\gamma_1 \gamma_1 = \gamma_2,~~~\gamma_2 \gamma_2 = \gamma_1,~~~
\gamma_1 \gamma_2 = \gamma_0,\\
&\rho_0 \rho_0 = \frac{1}{3}\gamma_0 + \frac{1}{3}\gamma_1 + \frac{1}{3}\gamma_2,~~~
\gamma_1 \rho_0 = \rho_0,~~~\gamma_2 \rho_0 = \rho_0.  
\end{align*}

\bigskip
\noindent
{\bf 4.3 } \ The case that $G$ is the symmetric group 
$ S_3 = \mathbb{Z}_3 \rtimes_\alpha \mathbb{Z}_2$ of degree 3 and  $G_0 = \mathbb{Z}_2$.
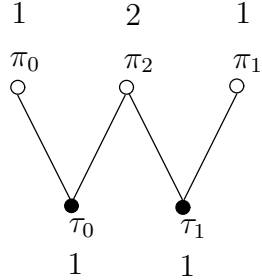
\begin{center}
{\unitlength 0.1in%
\begin{picture}(12.1200,13.2000)(19.8000,-19.6000)%
%
\special{pn 8}%
\special{pa 2022 1146}%
\special{pa 2300 1702}%
\special{fp}%
\special{pa 2307 1695}%
\special{pa 2578 1146}%
\special{fp}%
\put(19.8000,-10.2000){\makebox(0,0)[lb]{$\pi_0$}}%
\put(31.3500,-10.2700){\makebox(0,0)[lb]{$\pi_1$}}%
\put(22.8600,-18.8300){\makebox(0,0)[lb]{$\tau_0$}}%
%
\special{pn 8}%
\special{ar 2585 1104 36 36 0.0000000 6.2831853}%
%
\special{pn 8}%
\special{ar 2022 1104 36 36 0.0000000 6.2831853}%
%
\special{sh 1.000}%
\special{ia 2300 1723 35 35 0.0000000 6.2831853}%
\special{pn 8}%
\special{ar 2300 1723 35 35 0.0000000 6.2831853}%
%
\special{pn 8}%
\special{pa 2599 1139}%
\special{pa 2877 1696}%
\special{fp}%
\special{pa 2884 1689}%
\special{pa 3156 1139}%
\special{fp}%
%
\special{pn 8}%
\special{ar 3156 1097 36 36 0.0000000 6.2831853}%
%
\special{sh 1.000}%
\special{ia 2877 1731 36 36 0.0000000 6.2831853}%
\special{pn 8}%
\special{ar 2877 1731 36 36 0.0000000 6.2831853}%
\put(28.6300,-18.9000){\makebox(0,0)[lb]{$\tau_1$}}%
\put(25.7100,-10.2700){\makebox(0,0)[lb]{$\pi_2$}}%
\put(19.9000,-7.7000){\makebox(0,0)[lb]{1}}%
\put(31.5000,-7.7000){\makebox(0,0)[lb]{1}}%
\put(25.8000,-7.8000){\makebox(0,0)[lb]{2}}%
\put(22.8000,-20.8000){\makebox(0,0)[lb]{1}}%
\put(28.6000,-20.9000){\makebox(0,0)[lb]{1}}%
\end{picture}}%
\end{center}

\bigskip
\noindent
$\mathcal{K}(\hat{G} \cup \widehat{G_0}) 
= \{(ch(\pi_i), \circ), (ch(\tau_j), \bullet): \pi_i \in \hat{G}, \tau_j \in \widehat{G_0}\}$. 
Put $\gamma_i = (ch(\pi_i), \circ)$ and $\rho_j = (ch(\tau_j), \bullet)$.  
Then the structure equations are
\begin{align*}
&\gamma_1 \gamma_1 = \gamma_0,~~~
\gamma_2 \gamma_2 = \frac{1}{4}\gamma_0 + \frac{1}{4}\gamma_1 + \frac{1}{2}\gamma_2,~~~
\gamma_1 \gamma_2 = \gamma_2,~~~
\rho_0 \rho_0 = \rho_1 \rho_1 = \frac{1}{3}\gamma_0 + \frac{2}{3}\gamma_2,\\
&\rho_0 \rho_1 = \rho_1 \rho_0 = \frac{1}{3}\gamma_1 + \frac{2}{3}\gamma_2,~~~
\gamma_0 \rho_0 = \rho_0,~~~\gamma_1 \rho_0 = \rho_1,~~~
\gamma_2 \rho_0 = \frac{1}{2}\rho_0 + \frac{1}{2}\rho_1,\\
&\gamma_0 \rho_1 = \rho_1,~~~\gamma_1 \rho_1 = \rho_0,~~~
\gamma_2 \rho_1 = \frac{1}{2}\rho_0 + \frac{1}{2}\rho_1. 
\end{align*}

\medskip
\noindent
{\bf Remark } $\mathcal{K}(\widehat{\mathbb{Z}_2} \cup \widehat{\{e\}}) = 
\mathcal{K}(A_3)$, 
$\mathcal{K}(\widehat{\mathbb{Z}_3} \cup \widehat{\{e\}}) = 
\mathcal{K}(D_4)$ and 
$\mathcal{K}(\widehat{S_3} \cup \widehat{\mathbb{Z}_2}) = \mathcal{K}(A_5)$ 
where $\mathcal{K}(A_3)$, $\mathcal{K}(D_4)$ and 
$\mathcal{K}(A_5)$ are Sunder-Wildberger's hypergroups ([SW]) associated with  
Dynkin diagrams of type $A_3$, $D_4$ and $A_5$ respectively. 

\bigskip
\noindent
{\bf 4.4 } \ The case that 
$G = \mathbb{Z}_4 = \{e,g, g^2, g^3\}$ $(g^4 = e)$ and $G_0 = \mathbb{Z}_2 $.
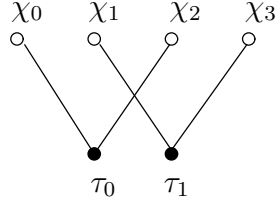
\begin{center}
{\unitlength 0.1in%
\begin{picture}(12.6400,9.1000)(5.6800,-18.9000)%
\put(5.7000,-11.1000){\makebox(0,0)[lb]{$\chi_0$}}%
\put(9.7000,-11.1000){\makebox(0,0)[lb]{$\chi_1$}}%
\put(9.8000,-20.2000){\makebox(0,0)[lb]{$\tau_0$}}%
\put(13.6000,-20.2000){\makebox(0,0)[lb]{$\tau_1$}}%
\put(13.9000,-11.1000){\makebox(0,0)[lb]{$\chi_2$}}%
%
\special{pn 8}%
\special{ar 1000 1200 32 32 0.0000000 6.2831853}%
%
\special{sh 1.000}%
\special{ia 1400 1800 32 32 0.0000000 6.2831853}%
\special{pn 8}%
\special{ar 1400 1800 32 32 0.0000000 6.2831853}%
%
\special{pn 8}%
\special{ar 1400 1200 32 32 0.0000000 6.2831853}%
%
\special{sh 1.000}%
\special{ia 1000 1800 32 32 0.0000000 6.2831853}%
\special{pn 8}%
\special{ar 1000 1800 32 32 0.0000000 6.2831853}%
%
\special{pn 8}%
\special{ar 1800 1200 32 32 0.0000000 6.2831853}%
%
\special{pn 8}%
\special{ar 600 1200 32 32 0.0000000 6.2831853}%
%
\special{pn 8}%
\special{pa 1790 1230}%
\special{pa 1400 1780}%
\special{fp}%
%
\special{pn 8}%
\special{pa 1390 1230}%
\special{pa 1000 1800}%
\special{fp}%
%
\special{pn 8}%
\special{pa 640 1230}%
\special{pa 1000 1800}%
\special{fp}%
%
\special{pn 8}%
\special{pa 1030 1230}%
\special{pa 1400 1780}%
\special{fp}%
\put(17.8000,-11.1000){\makebox(0,0)[lb]{$\chi_3$}}%
\end{picture}}%
\end{center}

\bigskip

$\mathcal{K}(\hat{G} \cup \widehat{G_0}) 
= \{(ch(\chi_i), \circ), (ch(\tau_j), \bullet): \chi_i \in \hat{G}, \tau_j \in \widehat{G_0}\}$. 
Put $\gamma_i = (ch(\chi_i), \circ)$ and $\rho_j = (ch(\tau_j), \bullet)$.  
Then the structure equations are
\begin{align*}
&\gamma_1 \gamma_1 = \gamma_2,~~~\gamma_2 \gamma_2 = \gamma_0,~~~
\gamma_3 \gamma_3 = \gamma_1,~~~
\gamma_1 \gamma_2 = \gamma_3,~~~
\gamma_1 \gamma_3 = \gamma_0,~~~
\gamma_2 \gamma_3 = \gamma_1,\\
&\rho_0 \rho_0 = \rho_1 \rho_1 = \frac{1}{2}\gamma_0 + \frac{1}{2}\gamma_2,~~~
\rho_0 \rho_1 = \rho_1 \rho_0 = \frac{1}{2}\gamma_1 + \frac{1}{2}\gamma_3,~~~
\gamma_0 \rho_0 = \rho_0,~~~\gamma_1 \rho_0 = \rho_1,\\
&\gamma_2 \rho_0 = \rho_0,~~~\gamma_3 \rho_0 = \rho_1,~~~
\gamma_0 \rho_1 = \rho_1,~~~\gamma_1 \rho_1 = \rho_0,~~~
\gamma_2 \rho_1 = \rho_1,~~~\gamma_3 \rho_1 = \rho_0. 
\end{align*}

\bigskip
\noindent
{\bf 4.5 } \ The case that 
$G = \mathbb{Z}_2 \times \mathbb{Z}_2 = \{(e,e),(e,g), (g,e), (g,g)\}$ $(g^2 = e)$ and $G_0 = \mathbb{Z}_2 $.
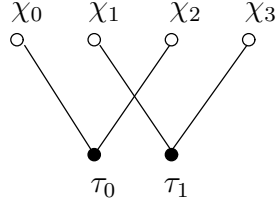
\begin{center}
{\unitlength 0.1in%
\begin{picture}(12.6400,9.1000)(5.6800,-18.9000)%
\put(5.7000,-11.1000){\makebox(0,0)[lb]{$\chi_0$}}%
\put(9.7000,-11.1000){\makebox(0,0)[lb]{$\chi_1$}}%
\put(9.8000,-20.2000){\makebox(0,0)[lb]{$\tau_0$}}%
\put(13.6000,-20.2000){\makebox(0,0)[lb]{$\tau_1$}}%
\put(13.9000,-11.1000){\makebox(0,0)[lb]{$\chi_2$}}%
%
\special{pn 8}%
\special{ar 1000 1200 32 32 0.0000000 6.2831853}%
%
\special{sh 1.000}%
\special{ia 1400 1800 32 32 0.0000000 6.2831853}%
\special{pn 8}%
\special{ar 1400 1800 32 32 0.0000000 6.2831853}%
%
\special{pn 8}%
\special{ar 1400 1200 32 32 0.0000000 6.2831853}%
%
\special{sh 1.000}%
\special{ia 1000 1800 32 32 0.0000000 6.2831853}%
\special{pn 8}%
\special{ar 1000 1800 32 32 0.0000000 6.2831853}%
%
\special{pn 8}%
\special{ar 1800 1200 32 32 0.0000000 6.2831853}%
%
\special{pn 8}%
\special{ar 600 1200 32 32 0.0000000 6.2831853}%
%
\special{pn 8}%
\special{pa 1790 1230}%
\special{pa 1400 1780}%
\special{fp}%
%
\special{pn 8}%
\special{pa 1390 1230}%
\special{pa 1000 1800}%
\special{fp}%
%
\special{pn 8}%
\special{pa 640 1230}%
\special{pa 1000 1800}%
\special{fp}%
%
\special{pn 8}%
\special{pa 1030 1230}%
\special{pa 1400 1780}%
\special{fp}%
\put(17.8000,-11.1000){\makebox(0,0)[lb]{$\chi_3$}}%
\end{picture}}%
\end{center}

\bigskip

$\mathcal{K}(\hat{G} \cup \widehat{G_0}) 
= \{(ch(\chi_i), \circ), (ch(\tau_j), \bullet): \chi_i \in \hat{G}, \tau_j \in \widehat{G_0}\}$. 
Put $\gamma_i = (ch(\chi_i), \circ)$ and $\rho_j = (ch(\tau_j), \bullet)$.  
Then the structure equations are
\begin{align*}
&\gamma_1 \gamma_1 = \gamma_0,~~~\gamma_2 \gamma_2 = \gamma_0,~~~
\gamma_3 \gamma_3 = \gamma_0,~~~
\gamma_1 \gamma_2 = \gamma_3,~~~
\gamma_1 \gamma_3 = \gamma_2,~~~
\gamma_2 \gamma_3 = \gamma_1,\\
&\rho_0 \rho_0 = \rho_1 \rho_1 = \frac{1}{2}\gamma_0 + \frac{1}{2}\gamma_2,~~~
\rho_0 \rho_1 = \rho_1 \rho_0 = \frac{1}{2}\gamma_1 + \frac{1}{2}\gamma_3,~~~
\gamma_0 \rho_0 = \rho_0,~~~\gamma_1 \rho_0 = \rho_1,\\
&\gamma_2 \rho_0 = \rho_0,~~~\gamma_3 \rho_0 = \rho_1,~~~
\gamma_0 \rho_1 = \rho_1,~~~\gamma_1 \rho_1 = \rho_0,~~~
\gamma_2 \rho_1 = \rho_1,~~~\gamma_3 \rho_1 = \rho_0. 
\end{align*}

\bigskip

\noindent
{\bf Remark} \ We note that Frobenius diagrams of 4.4 and 4.5 are same but their  
hypergroup structures are different. 

\bigskip
\noindent
{\bf 4.6 } \ The case that 
$G = S_3 = \mathbb{Z}_3 \rtimes_\alpha \mathbb{Z}_2$ and $G_0 = \mathbb{Z}_3$.
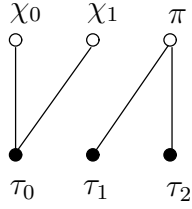
\begin{center}
{\unitlength 0.1in%
\begin{picture}(  8.6400,  9.3000)(  5.6800,-19.1000)%
\put(5.7000,-11.1000){\makebox(0,0)[lb]{$\chi_0$}}%
\put(9.7000,-11.1000){\makebox(0,0)[lb]{$\chi_1$}}%
\put(5.7000,-20.3000){\makebox(0,0)[lb]{$\tau_0$}}%
\put(9.5000,-20.3000){\makebox(0,0)[lb]{$\tau_1$}}%
\put(13.9000,-11.1000){\makebox(0,0)[lb]{$\pi$}}%
%
\special{pn 8}%
\special{ar 1000 1200 32 32  0.0000000  6.2831853}%
%
\special{sh 1.000}%
\special{ia 600 1800 32 32  0.0000000  6.2831853}%
\special{pn 8}%
\special{ar 600 1800 32 32  0.0000000  6.2831853}%
%
\special{pn 8}%
\special{ar 1400 1200 32 32  0.0000000  6.2831853}%
%
\special{sh 1.000}%
\special{ia 1000 1800 32 32  0.0000000  6.2831853}%
\special{pn 8}%
\special{ar 1000 1800 32 32  0.0000000  6.2831853}%
%
\special{sh 1.000}%
\special{ia 1400 1800 32 32  0.0000000  6.2831853}%
\special{pn 8}%
\special{ar 1400 1800 32 32  0.0000000  6.2831853}%
%
\special{pn 8}%
\special{ar 600 1200 32 32  0.0000000  6.2831853}%
%
\special{pn 8}%
\special{pa 600 1240}%
\special{pa 600 1790}%
\special{fp}%
\special{pa 600 1790}%
\special{pa 600 1770}%
\special{fp}%
%
\special{pn 8}%
\special{pa 990 1240}%
\special{pa 600 1790}%
\special{fp}%
%
\special{pn 8}%
\special{pa 1390 1230}%
\special{pa 1000 1800}%
\special{fp}%
%
\special{pn 8}%
\special{pa 1410 1240}%
\special{pa 1410 1780}%
\special{fp}%
\put(13.8000,-20.4000){\makebox(0,0)[lb]{$\tau_2$}}%
\end{picture}}%
\end{center}

\bigskip

\noindent
$\mathcal{K}(\hat{G} \cup \widehat{G_0})$ is not a hypergroup by Corollary 3.17.

\bigskip
\noindent
{\bf 4.7 } \ The case that $G$ is the dihedral group 
$D_4 = \mathbb{Z}_4 \rtimes_\alpha \mathbb{Z}_2$ and $G_0 = \mathbb{Z}_2$.
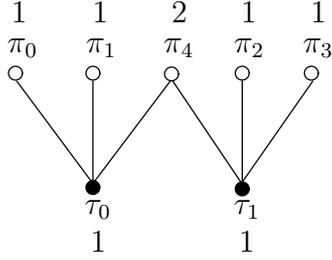
\begin{center}
{\unitlength 0.1in%
\begin{picture}(16.1000,12.0000)(19.8000,-19.1000)%
\put(19.8000,-10.2100){\makebox(0,0)[lb]{$\pi_0$}}%
\put(23.8800,-10.2100){\makebox(0,0)[lb]{$\pi_1$}}%
\put(23.8100,-18.6500){\makebox(0,0)[lb]{$\tau_0$}}%
%
\special{pn 8}%
\special{ar 2423 1106 36 36 0.0000000 6.2831853}%
%
\special{pn 8}%
\special{ar 2022 1106 36 36 0.0000000 6.2831853}%
%
\special{sh 1.000}%
\special{ia 2423 1703 36 36 0.0000000 6.2831853}%
\special{pn 8}%
\special{ar 2423 1703 36 36 0.0000000 6.2831853}%
%
\special{pn 8}%
\special{ar 2830 1110 36 36 0.0000000 6.2831853}%
%
\special{pn 8}%
\special{pa 2423 1148}%
\special{pa 2423 1682}%
\special{fp}%
\special{pa 2823 1148}%
\special{pa 2423 1682}%
\special{fp}%
\special{pa 2029 1148}%
\special{pa 2423 1682}%
\special{fp}%
\put(31.5400,-10.2100){\makebox(0,0)[lb]{$\pi_2$}}%
%
\special{pn 8}%
\special{ar 3196 1106 36 36 0.0000000 6.2831853}%
%
\special{pn 8}%
\special{ar 3554 1106 36 36 0.0000000 6.2831853}%
%
\special{sh 1.000}%
\special{ia 3196 1710 36 36 0.0000000 6.2831853}%
\special{pn 8}%
\special{ar 3196 1710 36 36 0.0000000 6.2831853}%
%
\special{pn 8}%
\special{pa 2837 1148}%
\special{pa 3189 1675}%
\special{fp}%
\special{pa 3196 1148}%
\special{pa 3196 1675}%
\special{fp}%
\special{pa 3561 1148}%
\special{pa 3196 1675}%
\special{fp}%
\put(27.8800,-10.2100){\makebox(0,0)[lb]{$\pi_4$}}%
\put(35.1200,-10.2100){\makebox(0,0)[lb]{$\pi_3$}}%
\put(31.5400,-18.6500){\makebox(0,0)[lb]{$\tau_1$}}%
\put(20.0000,-8.4000){\makebox(0,0)[lb]{1}}%
\put(24.2000,-8.4000){\makebox(0,0)[lb]{1}}%
\put(28.3000,-8.4000){\makebox(0,0)[lb]{2}}%
\put(31.9000,-8.4000){\makebox(0,0)[lb]{1}}%
\put(35.5000,-8.4000){\makebox(0,0)[lb]{1}}%
\put(24.1000,-20.4000){\makebox(0,0)[lb]{1}}%
\put(31.8000,-20.4000){\makebox(0,0)[lb]{1}}%
\end{picture}}%
\end{center}

\bigskip

$\mathcal{K}(\hat{G} \cup \widehat{G_0}) 
= \{(ch(\pi_i), \circ), (ch(\tau_j), \bullet): \pi_i \in \hat{G}, \tau_j \in \widehat{G_0}\}$. 
Put $\gamma_i = (ch(\pi_i), \circ)$ and $\rho_j = (ch(\tau_j), \bullet)$.  
Then the structure equations are
\begin{align*}
&\gamma_1 \gamma_1 = \gamma_0,~~~ 
\gamma_2 \gamma_2 
= \frac{1}{4}\gamma_0 + \frac{1}{4}\gamma_1 + \frac{1}{4}\gamma_3 + \frac{1}{4}\gamma_4,~~~
\gamma_3 \gamma_3 = \gamma_0,~~~
\gamma_4 \gamma_4 = \gamma_0,\\
&\gamma_1 \gamma_2 = \gamma_2,~~~
\gamma_1 \gamma_3 = \gamma_4,~~~
\gamma_1 \gamma_4 = \gamma_3,~~~
\gamma_2 \gamma_3 = \gamma_2,~~~
\gamma_2 \gamma_4 = \gamma_2,~~~
\gamma_3 \gamma_4 = \gamma_1,\\
&\rho_0 \rho_0 = \rho_1 \rho_1 = \frac{1}{4}\gamma_0 + \frac{1}{4}\gamma_1 + \frac{1}{2}\gamma_2,~~~
\rho_0 \rho_1 = \frac{1}{2}\gamma_2 + \frac{1}{4}\gamma_3 + \frac{1}{4}\gamma_4,\\
&\gamma_1 \rho_0 = \rho_0,~~~
\gamma_2 \rho_0 = \frac{1}{2}\rho_0 + \frac{1}{2}\rho_1,
\gamma_3 \rho_0 = \rho_1,~~~
\gamma_4 \rho_0 = \rho_1,\\
&\gamma_1 \rho_1 = \rho_1,~~~
\gamma_2 \rho_1 = \frac{1}{2}\rho_0 + \frac{1}{2}\rho_1,~~~
\gamma_3 \rho_1 = \rho_0,~~~
\gamma_4 \rho_1 = \rho_0.
\end{align*}

\bigskip
\noindent
{\bf 4.8 } \ The case that $G$ is the alternating group 
$A_4 = (\mathbb{Z}_2 \times \mathbb{Z}_2) \rtimes_\alpha \mathbb{Z}_3$ of degree 4 and  $G_0 = \mathbb{Z}_3$.
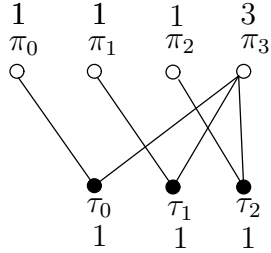
\begin{center}
{\unitlength 0.1in%
\begin{picture}(12.5200,11.7000)(19.8000,-19.0000)%
\put(19.8000,-10.2100){\makebox(0,0)[lb]{$\pi_0$}}%
\put(23.8800,-10.2100){\makebox(0,0)[lb]{$\pi_1$}}%
\put(23.8100,-18.6500){\makebox(0,0)[lb]{$\tau_0$}}%
%
\special{pn 8}%
\special{ar 2423 1106 36 36 0.0000000 6.2831853}%
%
\special{pn 8}%
\special{ar 2022 1106 36 36 0.0000000 6.2831853}%
%
\special{sh 1.000}%
\special{ia 2423 1703 36 36 0.0000000 6.2831853}%
\special{pn 8}%
\special{ar 2423 1703 36 36 0.0000000 6.2831853}%
%
\special{pn 8}%
\special{ar 2830 1110 36 36 0.0000000 6.2831853}%
\put(27.9000,-10.1000){\makebox(0,0)[lb]{$\pi_2$}}%
%
\special{pn 8}%
\special{ar 3196 1106 36 36 0.0000000 6.2831853}%
%
\special{sh 1.000}%
\special{ia 2830 1710 36 36 0.0000000 6.2831853}%
\special{pn 8}%
\special{ar 2830 1710 36 36 0.0000000 6.2831853}%
%
\special{sh 1.000}%
\special{ia 3196 1710 36 36 0.0000000 6.2831853}%
\special{pn 8}%
\special{ar 3196 1710 36 36 0.0000000 6.2831853}%
\put(31.8000,-10.1000){\makebox(0,0)[lb]{$\pi_3$}}%
\put(31.5400,-18.6500){\makebox(0,0)[lb]{$\tau_2$}}%
\put(28.0000,-18.8000){\makebox(0,0)[lb]{$\tau_1$}}%
%
\special{pn 8}%
\special{pa 2030 1140}%
\special{pa 2423 1703}%
\special{fp}%
%
\special{pn 8}%
\special{pa 2440 1140}%
\special{pa 2830 1710}%
\special{fp}%
%
\special{pn 8}%
\special{pa 2860 1150}%
\special{pa 3196 1710}%
\special{fp}%
%
\special{pn 8}%
\special{pa 3170 1130}%
\special{pa 2423 1703}%
\special{fp}%
%
\special{pn 8}%
\special{pa 3170 1130}%
\special{pa 2830 1710}%
\special{fp}%
%
\special{pn 8}%
\special{pa 3170 1130}%
\special{pa 3196 1710}%
\special{fp}%
\put(19.9000,-8.6000){\makebox(0,0)[lb]{1}}%
\put(19.9000,-8.6000){\makebox(0,0)[lb]{1}}%
\put(24.1000,-8.6000){\makebox(0,0)[lb]{1}}%
\put(24.1000,-8.6000){\makebox(0,0)[lb]{1}}%
\put(28.1000,-8.7000){\makebox(0,0)[lb]{1}}%
\put(24.1000,-20.2000){\makebox(0,0)[lb]{1}}%
\put(28.3000,-20.3000){\makebox(0,0)[lb]{1}}%
\put(31.8000,-20.3000){\makebox(0,0)[lb]{1}}%
\put(31.8000,-8.6000){\makebox(0,0)[lb]{3}}%
\end{picture}}%
\end{center}
$\mathcal{K}(\hat{G} \cup \widehat{G_0}) 
= \{(ch(\pi_i), \circ), (ch(\tau_j), \bullet): \pi_i \in \hat{G}, \tau_j \in \widehat{G_0}\}$. 
Put $\gamma_i = (ch(\pi_i), \circ)$ and $\rho_j = (ch(\tau_j), \bullet)$.
Then the structure equations are
\begin{align*}
&\gamma_1 \gamma_1 = \gamma_2,~~~
\gamma_2 \gamma_2 = \gamma_1,~~~
\gamma_3 \gamma_3 
= \frac{1}{9}\gamma_0 + \frac{1}{9}\gamma_1 + \frac{1}{9}\gamma_2 + \frac{2}{3}\gamma_3,~~~
\gamma_1 \gamma_2 = \gamma_0,~~~
\gamma_1 \gamma_3 = \gamma_3,\\
&\gamma_2 \gamma_3 = \gamma_3,~~~
\rho_0 \rho_0 = \rho_1 \rho_2 = \frac{1}{4}\gamma_0 +  \frac{3}{4}\gamma_3,~~~
\rho_0 \rho_1 = \frac{1}{4}\gamma_1 + \frac{3}{4}\gamma_3,~~~
\rho_0 \rho_2 = \frac{1}{4}\gamma_2 + \frac{3}{4}\gamma_3,\\
&\gamma_1 \rho_0 = \rho_1,~~~
\gamma_2 \rho_0 = \rho_2,~~~
\gamma_1 \rho_1 = \rho_2,~~~
\gamma_2 \rho_1 = \rho_0, \\
&\gamma_3 \rho_0 = \gamma_3 \rho_1 = \gamma_3 \rho_2 = \frac{1}{3}\rho_0 + \frac{1}{3}\rho_1 + \frac{1}{3}\rho_2. 
\end{align*}

\bigskip
\noindent
{\bf 4.9 } \ The case that $G$ is the symmetric group 
$S_4 = A_4 \rtimes_\alpha \mathbb{Z}_2$ of degree 4 and $G_0 = \mathbb{Z}_2$.
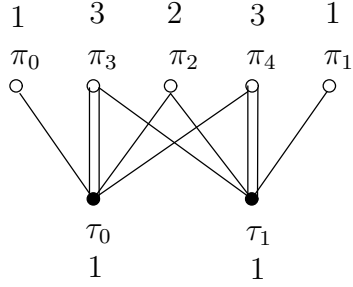
\begin{center}
{\unitlength 0.1in%
\begin{picture}(16.8600,13.4000)(5.6800,-20.8000)%
\put(5.7000,-11.1000){\makebox(0,0)[lb]{$\pi_0$}}%
\put(9.7000,-11.1000){\makebox(0,0)[lb]{$\pi_3$}}%
\put(9.6000,-20.2000){\makebox(0,0)[lb]{$\tau_0$}}%
\put(17.9000,-20.3000){\makebox(0,0)[lb]{$\tau_1$}}%
\put(13.9000,-11.1000){\makebox(0,0)[lb]{$\pi_2$}}%
%
\special{pn 8}%
\special{ar 1000 1200 32 32 0.0000000 6.2831853}%
%
\special{sh 1.000}%
\special{ia 1000 1790 32 32 0.0000000 6.2831853}%
\special{pn 8}%
\special{ar 1000 1790 32 32 0.0000000 6.2831853}%
%
\special{pn 8}%
\special{ar 1400 1200 32 32 0.0000000 6.2831853}%
%
\special{sh 1.000}%
\special{ia 1820 1790 32 32 0.0000000 6.2831853}%
\special{pn 8}%
\special{ar 1820 1790 32 32 0.0000000 6.2831853}%
%
\special{pn 8}%
\special{ar 600 1200 32 32 0.0000000 6.2831853}%
\put(17.9200,-11.1000){\makebox(0,0)[lb]{$\pi_4$}}%
\put(21.9200,-11.1000){\makebox(0,0)[lb]{$\pi_1$}}%
%
\special{pn 8}%
\special{ar 2222 1200 32 32 0.0000000 6.2831853}%
%
\special{pn 8}%
\special{ar 1822 1200 32 32 0.0000000 6.2831853}%
%
\special{pn 8}%
\special{pa 1030 1210}%
\special{pa 1030 1800}%
\special{fp}%
%
\special{pn 8}%
\special{pa 1800 1200}%
\special{pa 1800 1800}%
\special{fp}%
%
\special{pn 8}%
\special{pa 1852 1210}%
\special{pa 1850 1800}%
\special{fp}%
%
\special{pn 8}%
\special{pa 980 1210}%
\special{pa 980 1780}%
\special{fp}%
%
\special{pn 8}%
\special{pa 1400 1240}%
\special{pa 1000 1790}%
\special{fp}%
%
\special{pn 8}%
\special{pa 1800 1230}%
\special{pa 1000 1790}%
\special{fp}%
%
\special{pn 8}%
\special{pa 1030 1210}%
\special{pa 1820 1790}%
\special{fp}%
%
\special{pn 8}%
\special{pa 1400 1240}%
\special{pa 1820 1790}%
\special{fp}%
%
\special{pn 8}%
\special{pa 2210 1230}%
\special{pa 1820 1790}%
\special{fp}%
%
\special{pn 8}%
\special{pa 620 1240}%
\special{pa 1000 1790}%
\special{fp}%
\put(5.7000,-8.9000){\makebox(0,0)[lb]{1}}%
\put(9.8000,-8.7000){\makebox(0,0)[lb]{3}}%
\put(13.8000,-8.7000){\makebox(0,0)[lb]{2}}%
\put(18.1000,-8.8000){\makebox(0,0)[lb]{3}}%
\put(22.0000,-8.7000){\makebox(0,0)[lb]{1}}%
\put(9.7000,-22.0000){\makebox(0,0)[lb]{1}}%
\put(18.1000,-22.1000){\makebox(0,0)[lb]{1}}%
\end{picture}}%
\end{center}

\bigskip
\medskip

$\mathcal{K}(\hat{G} \cup \widehat{G_0}) 
= \{(ch(\pi_i), \circ), (ch(\tau_j), \bullet): \pi_i \in \hat{G}, \tau_j \in \widehat{G_0}\}$. 
Put $\gamma_i = (ch(\pi_i), \circ)$ and $\rho_j = (ch(\tau_j), \bullet)$.  
Then the structure equations are
\begin{align*}
&\gamma_1  \gamma_1 = \gamma_0,~~~
\gamma_2  \gamma_2 = \frac{1}{4}\gamma_0 + \frac{1}{4}\gamma_1 + \frac{1}{2}\gamma_2,~~~
\gamma_3 \gamma_3 = \gamma_4  \gamma_4 
= \frac{1}{9}\gamma_0 + \frac{2}{9}\gamma_2 + \frac{1}{3}\gamma_3 + \frac{1}{3}\gamma_4,\\
&\gamma_1  \gamma_2 = \gamma_2,~~~\gamma_1  \gamma_3 = \gamma_4,~~~
\gamma_1 \gamma_4 = \gamma_3,~~~
\gamma_2 \gamma_3 = \gamma_2  \gamma_4 
= \frac{1}{2}\gamma_3 + \frac{1}{2}\gamma_4,\\
&\gamma_3 \gamma_4 
= \frac{1}{9}\gamma_1 + \frac{2}{9}\gamma_2 + \frac{1}{3}\gamma_3 + \frac{1}{3}\gamma_4,~~~
\rho_0 \rho_0 = \rho_1 \rho_1 = 
\frac{1}{12}\gamma_0 + \frac{1}{6}\gamma_2 + \frac{1}{2}\gamma_3 + \frac{1}{4}\gamma_4,\\
&\rho_0 \rho_1 = \rho_1 \rho_0 = 
\frac{1}{12}\gamma_1 + \frac{1}{6}\gamma_2 + \frac{1}{4}\gamma_3 + \frac{1}{2}\gamma_4,~~~
\gamma_0 \rho_0 = \rho_0,~~~\gamma_1 \rho_0 = \rho_1,\\
&\gamma_2 \rho_0 = \frac{1}{2}\rho_0 + \frac{1}{2}\rho_1,~~~  
\gamma_3 \rho_0 =  \frac{2}{3}\rho_0 + \frac{1}{3}\rho_1,~~~
\gamma_4 \rho_0 =  \frac{1}{3}\rho_0 + \frac{2}{3}\rho_1,~~~
\gamma_0 \rho_1 = \rho_1,\\
&\gamma_1 \rho_1 = \rho_0,~~~
\gamma_2 \rho_1 = \frac{1}{2}\rho_0 + \frac{1}{2}\rho_1,~~~  
\gamma_3 \rho_1 =  \frac{1}{3}\rho_0 + \frac{2}{3}\rho_1,~~~
\gamma_4 \rho_1 =  \frac{2}{3}\rho_0 + \frac{1}{3}\rho_1.          
\end{align*}

\bigskip
\noindent
{\bf 4.10 } \ The case that $G$ is the symmetric group 
$S_4$ of degree 4 and $G_0$ is the symmetric group $S_3$ of degree 3.
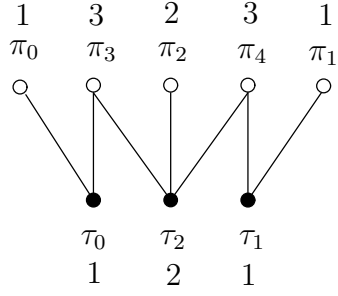
\begin{center}
{\unitlength 0.1in%
\begin{picture}(16.6400,13.7000)(17.7000,-19.1000)%
\put(17.7000,-8.6000){\makebox(0,0)[lb]{$\pi_0$}}%
\put(29.5000,-8.8000){\makebox(0,0)[lb]{$\pi_4$}}%
\put(21.4000,-18.5000){\makebox(0,0)[lb]{$\tau_0$}}%
%
\special{pn 8}%
\special{ar 2606 996 36 36 0.0000000 6.2831853}%
%
\special{pn 8}%
\special{ar 3398 1000 36 36 0.0000000 6.2831853}%
%
\special{sh 1.000}%
\special{ia 2606 1596 35 35 0.0000000 6.2831853}%
\special{pn 8}%
\special{ar 2606 1596 35 35 0.0000000 6.2831853}%
%
\special{sh 1.000}%
\special{ia 2206 1596 36 36 0.0000000 6.2831853}%
\special{pn 8}%
\special{ar 2206 1596 36 36 0.0000000 6.2831853}%
\put(25.5000,-18.5000){\makebox(0,0)[lb]{$\tau_2$}}%
\put(21.7000,-8.7000){\makebox(0,0)[lb]{$\pi_3$}}%
\put(33.2000,-8.9000){\makebox(0,0)[lb]{$\pi_1$}}%
%
\special{pn 8}%
\special{ar 2206 996 36 36 0.0000000 6.2831853}%
%
\special{pn 8}%
\special{ar 3006 996 36 36 0.0000000 6.2831853}%
%
\special{sh 1.000}%
\special{ia 3006 1596 36 36 0.0000000 6.2831853}%
\special{pn 8}%
\special{ar 3006 1596 36 36 0.0000000 6.2831853}%
\put(29.6000,-18.5000){\makebox(0,0)[lb]{$\tau_1$}}%
\put(25.4000,-8.7000){\makebox(0,0)[lb]{$\pi_2$}}%
%
\special{pn 8}%
\special{ar 1826 1006 36 36 0.0000000 6.2831853}%
%
\special{pn 8}%
\special{pa 1856 1046}%
\special{pa 2206 1596}%
\special{fp}%
%
\special{pn 8}%
\special{pa 2206 1036}%
\special{pa 2206 1596}%
\special{fp}%
%
\special{pn 8}%
\special{pa 2206 1036}%
\special{pa 2606 1596}%
\special{fp}%
%
\special{pn 8}%
\special{pa 2606 1036}%
\special{pa 2606 1596}%
\special{fp}%
%
\special{pn 8}%
\special{pa 3006 1036}%
\special{pa 3006 1596}%
\special{fp}%
%
\special{pn 8}%
\special{pa 3386 1036}%
\special{pa 3006 1596}%
\special{fp}%
%
\special{pn 8}%
\special{pa 3006 1036}%
\special{pa 2606 1596}%
\special{fp}%
\put(18.0000,-6.8000){\makebox(0,0)[lb]{1}}%
\put(21.8000,-6.8000){\makebox(0,0)[lb]{3}}%
\put(25.7000,-6.7000){\makebox(0,0)[lb]{2}}%
\put(29.8000,-6.7000){\makebox(0,0)[lb]{3}}%
\put(33.7000,-6.8000){\makebox(0,0)[lb]{1}}%
\put(21.7000,-20.3000){\makebox(0,0)[lb]{1}}%
\put(25.8000,-20.4000){\makebox(0,0)[lb]{2}}%
\put(29.7000,-20.4000){\makebox(0,0)[lb]{1}}%
\end{picture}}%
\end{center}

\bigskip

$\mathcal{K}(\hat{G} \cup \widehat{G_0}) 
= \{(ch(\pi_i), \circ), (ch(\tau_j), \bullet): \pi_i \in \hat{G}, \tau_j \in \widehat{G_0}\}$. 
Put $\gamma_i = (ch(\pi_i), \circ)$ and $\rho_j = (ch(\tau_j), \bullet)$.  
Then the structure equations are
\begin{align*}
&\gamma_1  \gamma_1 = \gamma_0,~~~
\gamma_2  \gamma_2 = \frac{1}{4}\gamma_0 + \frac{1}{4}\gamma_1 + \frac{1}{2}\gamma_2,~~~
\gamma_3 \gamma_3 = \gamma_4  \gamma_4 
= \frac{1}{9}\gamma_0 + \frac{2}{9}\gamma_2 + \frac{1}{3}\gamma_3 + \frac{1}{3}\gamma_4,\\
&\gamma_1  \gamma_2 = \gamma_2,~~~\gamma_1  \gamma_3 = \gamma_4,~~~
\gamma_1 \gamma_4 = \gamma_3,~~~
\gamma_2 \gamma_3 = \gamma_2  \gamma_4 
= \frac{1}{2}\gamma_3 + \frac{1}{2}\gamma_4,\\
&\gamma_3 \gamma_4 
= \frac{1}{9}\gamma_1 + \frac{2}{9}\gamma_2 + \frac{1}{3}\gamma_3 + \frac{1}{3}\gamma_4,~~~
\rho_0 \rho_0 = \rho_1 \rho_1 = \frac{1}{4}\gamma_0 + \frac{3}{4}\gamma_3,\\
& \rho_2 \rho_2 = \frac{1}{16}\gamma_0 + \frac{1}{16}\gamma_1 + \frac{1}{8}\gamma_2 + 
\frac{3}{8}\gamma_3 + \frac{3}{8}\gamma_4, ~~~
\rho_1 \rho_2 = \frac{1}{4}\gamma_2 + \frac{3}{8}\gamma_3 + \frac{4}{8}\gamma_4,\\
&\gamma_0 \rho_0 = \rho_0,~~~\gamma_1 \rho_0 = \rho_1,~~~\gamma_2 \rho_0 = \rho_2,~~~
\gamma_3 \rho_0 = \frac{1}{3}\rho_0 + \frac{2}{3}\rho_2,~~~
\gamma_4 \rho_0 = \frac{1}{3}\rho_1 + \frac{2}{3}\rho_2,\\
&\gamma_0 \rho_1 = \rho_1,~~~\gamma_1 \rho_1 = \rho_0,~~~\gamma_2 \rho_1 = \rho_2,~~~
\gamma_3 \rho_1 = \frac{1}{3}\rho_1 + \frac{2}{3}\rho_2,~~~
\gamma_4 \rho_1 = \frac{1}{3}\rho_0 + \frac{2}{3}\rho_2,\\
&\gamma_0 \rho_2 = \rho_2,~~~\gamma_1 \rho_2 = \rho_2,~~~
\gamma_2 \rho_2 = \frac{1}{4}\rho_0 + \frac{1}{4}\rho_1 + \frac{1}{2}\rho_2,\\
&\gamma_3 \rho_2 = \gamma_4 \rho_0 = \frac{1}{6}\rho_0 + \frac{1}{6}\rho_1 + \frac{2}{3}\rho_2. 
\end{align*}

\bigskip

\textbf{Addresses }

\medskip

Herbert Heyer : Universit\"{a}t T\"{u}bingen

Mathematisches Institut

Auf der Morgenstelle 10

72076, T\"{u}bingen

Germany

\medskip

e-mail : herbert.heyer@uni-tuebingen.de

\bigskip

Satoshi Kawakami : Nara University of Education

Department of Mathematics

Takabatake-cho

Nara, 630-8528

Japan

\medskip

e-mail : kawakami@nara-edu.ac.jp

\bigskip
Tatsuya Tsurii : Osaka Prefecture University 

Graduate School of Science

1-1 Gakuen-cho, Nakaku, Sakai 

Osaka, 599-8531

Japan

\medskip

e-mail : dw301003@edu.osakafu-u.ac.jp

\bigskip

Satoe Yamanaka : Nara Women's University 

Faculty of Science 

Kitauoya-higashimachi, 

Nara, 630-8506

Japan 

\medskip

e-mail : s.yamanaka516@gmail.com


\bigskip



\begin{thebibliography}{00000}



\bibitem[BH]{BH}
W. Bloom and H. Heyer : 
\textit{Harmonic Analysis of Probability Measures on Hypergroups}, 
Walter de Gruyter, de Gruyter Studies in Mathematics 20(1995).




\bibitem[H]
{H}
T. Hirai :
\textit{Classical method of constructing a complete set of
irreducible representations of semidirect product of a 
compact group with a finite group}, 
Prob. Math. Statistic Vol. 33, Fasc. 2 (2013),
353--362.

\bibitem[F]
{F}
G. B. Folland : 
\textit{A Course in Abstract Harmonic Analysis}, 
Textbooks in Mathematics (1994). 



\bibitem[GHJ]
{GHJ}
F. M. Goodman, P. de la Harpe and V. F. R. Jones : 
\textit{Coxeter Graphs and Towers of Algebras}, 
Springer-Verlag (1989).  

\bibitem[HKKK]
{HKKK}
H. Heyer, Y. Katayama, S. Kawakami  and K. Kawasaki: 
\textit{Extensions of finite commutative hypergroups}, 
Sci. Math. Jpn., e-2007, 127-139. 


\bibitem[HK1]
{HK3}
H. Heyer and S. Kawakami : 
\textit{A cohomology approach to the extension problem for commutative hypergroups}, 
Semigroup Forum Vol. 83(2011), 371-394. 



\bibitem[HK2]
{HK1}
H. Heyer and S. Kawakami : 
\textit{An imprimitivity theorem for representations 
of a semi-direct product hypergroup}, Journal of Lie Theory Vol. 24
(2014), 159--178. 

\bibitem[HK3]
{HK2}
H. Heyer and S. Kawakami : 
\textit{Hypergroup structures arising from certain dual objects of 
a hypergroup}, to appear in Journal of the Math. Soc. of Japan. 




\bibitem[HKTY1]
{HKTY1}
H. Heyer, S. Kawakami, T. Tsurii and S. Yamanaka: 
\textit{A commutative hypergroup associated with a hyperfield}, 
submitted to Archiv der Mathematic.

\bibitem[HKTY2]
{HKTY2}
H. Heyer, S. Kawakami, T. Tsurii and S. Yamanaka: 
\textit{Hypergroups related to a pair of compact hypergroups}, 
preprint.




\bibitem[HKY]
{HKY}
H. Heyer, S. Kawakami and S. Yamanaka : 
\textit{Characters of induced representations of a  compact hypergroup}, 
Monatsh. Math., 2015, online published.






\bibitem[J]
{J}
 R. Jewett : 
\textit{Spaces with an abstract convolution of measures},
Adv. Math., 18(1975), 1--101.


\bibitem[KTY]
{KTY} 
S. Kawakami,  T. Tsurii and S. Yamanaka : 
\textit{Deformations of finite hypergroups }, 
Sci. Math. Jpn., 2015, online published.


\bibitem[SW]
{SW}
V. S. Sunder and N. J. Wildberger : 
\textit{Fusion rule algebras and walks on graphs}, 
The Proceedings of the Fifth Ramanujan Symposium on Harmonic Analysis, Ed. K. Parthasarathy, Publ. of the Ramanujan Inst., No.6, (1999), 53-80. 



\end{thebibliography}
\end {document}